\documentclass[12pt,leqno]{article}
\usepackage{amsmath,amsfonts}
\usepackage[dvips]{graphics,color}
\usepackage{latexsym}
\usepackage{amssymb}
\usepackage{graphicx}
\usepackage{tikz}
\date{}
\begin{document}
\title{On area-minimizing Pfaffian varieties}
\maketitle
\thanks{}
\author{\centerline {Hongbin Cui, Xiaoxiang Jiao and Xiaowei Xu\footnote{The corresponding author}}}
\maketitle

\section *\small{{\centerline{Abstract}}
\medskip\noindent

\small{There are two significant families of minimal real matrix varieties: determinantal varieties and skew-symmetric determinantal varieties, the later ones are also known as Pfaffian varieties. In 1999, Kerckhove and Lawlor [Duke Math.J. 96(2),401--424,1999] proved that determinantal varieties are area-minimizing except for two families. In this paper we prove that all Pfaffian varieties are area-minimizing with the exception of Pfaffian hypersurfaces.



\tableofcontents

\section{Introduction}

\medskip\noindent

  
In geometric measure theory (GMT), area-minimizing surfaces in $\textbf{R}^{n}$ are integral currents (cf. \cite{FF60}) that globally minimize the area functional with a given boundary, i.e. they are solutions of generalized Plateau's problem (see \cite{Fe69}, \cite{M16}). For an area-minimizing  surface $S$, at every interior singular point $p$ there is an oriented tangent cone as the blowing-up limit of $S$ (cf. 4.3.16 and 5.4.3(6) in \cite{Fe69}), this cone is itself area-minimizing (cf. 5.4.3(5) in \cite{Fe69}). Therefore, for better understanding what kinds of singularities can occur in area-minimizing surfaces, it is important to know which cones can be area-minimizing. Meanwhile, many new examples of area-minimizing surfaces can be produced by perturbations of area-minimizing cones (cf. \cite{HS85}).



\medskip\noindent

To be precise, an area-minimizing cone $C$ is a locally rectifiable, homothetically invariant area-minimizing surface that starts at a fixed point (often the origin $O$) in $\textbf{R}^n$. It minimizes the area among all integral currents with the same boundary of its link $C\cap \textbf{S}_{O}^{n-1}(1)$, where $\textbf{S}_{O}^{n-1}(1)$ denote the unit sphere of $\textbf{R}^{n}$ that centered at $O$. Note that $O$ becomes a natural singularity of $C$ when it is nonplanar. For a submanifold $M$ of $\textbf{S}_{O}^{n-1}(1)$, the cone $C(M)$ over $M$ is the union of rays from $O$ that pass through points in $M$. 

\medskip\noindent

 
 \medskip\noindent

\textbf{1.1 Regular cones and non-regular cones.}   Area-minimizing cones with the simplest singularity (an isolated singularity at origin $O$) can be spanned by compact smooth, minimal submanifolds in unit spheres. In these cases, the singular set of $C$ is $\{O\}$, which are called \textit{regular cones}
by Hardt and Simon (\cite{HS85}). The links of regular cones (as described above), are compact smooth submanifolds in unit spheres. 

\medskip\noindent

The study of area-minimizing regular cones has a long history since the famous \textit{Simons cones} ($p=q=3$) and \textit{Lawson cones}
$C_{p,q}:=\{(x_{0},\ldots,x_{p+q+1})\;|\;q(x_{0}^2+\cdots+x_{p}^2)-p(x_{p+1}^2+\cdots+x_{p+q+1}^2)=0 \}$ which play important roles in the research of Bernstein problem (cf. \cite{Sim68}, \cite{BDG69}, etc}). They are cones over product of spheres, it is known that $C_{p,q}$ are area-minimizing if and only if $p+q\geq 7$ or $p+q=6,|p-q|\leq 3$ (\cite{Sim68}, \cite{BDG69}, \cite{L72}, \cite{Si74}, \cite{La91}), alternative proofs can be seen in \cite{Da04}, \cite{DP09}. Additionally, the area-minimizations for cones over arbitrary product of spheres was completely classified in \cite{La91} by Lawlor, where he provided a general method, called the \textit{curvature criterion}, to prove a regular cone's area-minimization. 
By using Lawlor's curvature criterion, many new area-minimizing regular cones are found, see \cite{La91}, \cite{Ke94}, \cite{Ka02}, \cite{XYZ18}, \cite{TZ20}, \cite{JC22}, \cite{JCX22}, etc. For more historical notes and methods of proof concerning on a given regular minimal cone's area-minimization, we refer the reader to Section 2 in \cite{La98}.

\medskip\noindent
\medskip\noindent

Very little is known for non-regular area-minimizing cones before. There are several classical examples. The celebrated \textit{angle conjecture} (\cite{M83}) states that the union of two oriented $k$-planes are area-minimizing if they satisfy a condition about characterizing angles. This condition was confirmed to be both sufficient and necessary by Nance, Lawlor (\cite{Ma87}, \cite{La89}). These minimizing cones have self-intersected like singularities at origins. There are plenty of non-regular area-minimizing cones naturally arising from singular complex projective varieties in algebraic geometry, when considered as cones in ambient Euclidean spaces. Via calibrated geometry theory, complex algebraic varieties in  ambient spaces $\textbf{C}^{n}\cong \textbf{R}^{2n}$ are area-minimizing in any of their compact portions (\cite{HL82}, \cite{M16}). Often these singular varieties carry stratified singular sets (cf. \cite{W57}, \cite{W65-1}, \cite{W65-2}).

\medskip\noindent

It is not easy to detect a given real algebraic variety's area-minimization, even it is minimal. Recently, the minimality of some real matrix varieties in ambient Euclidean spaces attracts people's attention. Let $M(p,q;\textbf{R})$, $(p\leq q)$, be the space of $p\times q$ matrices over $\textbf{R}$, and let $Skew(n,\textbf{R})$ be the space of $n \times n$ skew-symmetric matrices over $\textbf{R}$. Define \textit{determinantal variety} $C(p,q,r)$, $(r< p)$,  (cf. \cite{BV06}, \cite{Ha92}) as

$$
C(p,q,r)=\{X \in M(p,q;\textbf{R}): {\rm rank} X \leq r \},
$$
its skew-symmetric analogs ${\bf C}(n,2r)$, $(2r<n)$, also known as \textit{Pfaffian variety}, are defined as
$$
{\bf C}(n,2r)=\{X \in Skew(n,\textbf{R}): {\rm rank} X \leq 2r \}.
$$
By studying Hsiang's minimal surface equation (cf. \cite{Hs67}) for algebraic hypersurfaces, $C(p,p,p-1)$, ${\bf C}(2n,2n-2)$ are proved to be minimal in \cite{T10-1}, \cite{HLT17}, respectively. For the general $C(p,q,r)$ and ${\bf C}(n,2r)$, their minimalities are proved by various methods in \cite{BCH21} \cite{Ko21}, \cite{Ko21}, respectively. 
 
\medskip\noindent
 
In \cite{KL99}, Kerckhove and Lawlor proved that those $C(p,q,r)$, $(p+q-2r\geq 4)$, are stratified area-minimizing cones. To overcome the barriers of large singular sets\footnote{Rank one cases: $C(p,q,1)$, also known as Segre varieties, are regular cones, the area-minimizations had been studied in \cite{Pe93} by using Lawlor's Curvature Criterion}, they used the area-coarea formula and slicing techniques from GMT (\cite{Fe69}, \cite{M16}), by showing each slices minimized weighted measure in associated slicings. This method now is known as \textit{directed slicing} by Lawlor in \cite{La98}. By using this method, Lawlor gave a new and simple proof for the area-minimization of $C(S^{3}\times S^{3})$ (see 3.4 in \cite{La98}); Lawlor and Morgan (\cite{LM96}) proved that three minimal surfaces meeting at $120^{\circ}$ minimize the area locally; and others can be found in \cite{La98}, \cite{LM96}. 
  
\medskip\noindent
\medskip\noindent



\textbf{1.2 Statements of main results.} Inspiring by the works \cite{La98} and \cite{KL99}, for the Pfaffian varieties, we prove 

\medskip\noindent
 
 \textbf{Main Theorem}. {\it The Pfaffian varieties are stratified area-minimizing cones except for one family of hypersurfaces}.

\medskip\noindent

Our main theorem gives new area-minimizing non-regular cones (algebraic varieties) which may have interesting geometric meanings (see Remark 1.2 below). We also give another proof for the minimality of Pfaffian varieties. 

\medskip\noindent

The proof of our main theorem relies on a symmetry reduction that $SO(n)$ has a canonical adjoint action on ${\bf C}(n,2r)$. Those distinct rotated images of bounded (by the focal radii) normal spaces of ${\bf C}(n,2r)\subset Skew(n,\textbf{R})$ at regular points are the \textit{primary slicings} as we needed. The intersections of those slicings with ${\bf C}(n,2r)$ are \textit{primary slices}, and every comparison surface with the truncated part of ${\bf C}(n,2r)$ will intersect with primary slicings almost everywhere. We continue slice up primary slicings by parabolic curves, which is called the \textit{secondary slicings}. By computing primary weighting functions $w_{1}$ and composite weighting functions $w_{1}w_{2}$, we show that each slice of ${\bf C}(n,2r)$ minimizes the weighted measure in associated slicing set, then the whole minimization result can be concluded by area-coarea formula.

\medskip\noindent

\textbf{Remark 1.1}.  In \cite{KL99}, the area-minimizing determinantal varieties $C(p,q,r)$, $(p+q-2r\geq 4)$, are in fact divided into two cases: (1) when $ p+q $ is even, $C(p,q,r)$ is orientable, then $ C(p,q,r) $ is area-minimizing; (2) when $ p+q $ is odd, $C(p,q,r)$ is unorientable, then $C(p,q,r)$ is area-minimizing in the sense of modulo 2 (see Proposition 4.6 and Theorem 8.3 in \cite{KL99}). However, Pfaffian varieties are all orientable (see Proposition 2.3), therefore there has no mod 2 area-minimizing cases in our main thereom. 

\medskip\noindent
\medskip\noindent

\textbf{1.3 Set-theoretical tangent cones of Pfaffian varieties.} Since Lawlor's directed slicing method only works for those ${\bf C}(n,2r)$, $(n-2r\geq 3)$, (see Proposition 5.2), we can not show the minimization of the Pfaffian hypersurfaces  ${\bf C}(2n,2n-2)$. However, by the set-theoretical tangent cone (cf. 3.1.21 in \cite{Fe69}) analysis, we can still talk about some properties of ${\bf C}(2n,2n-2)$. 

\medskip\noindent

{\bf Definition} (Tangent cone, see 3.1.21 in \cite{Fe69}).\,\,  For an normed vector space $X$, a subset $S\subset X$ and a point $a\in X$, the tangent cone of $S$ at $a$ is defined by
\begin{equation}\notag
{\rm Tan}(S,a)=\{v\in X: \forall \varepsilon >0, \exists x \in S  \ and \ r\in \textbf{R}^{+}, s.t. \ |x-a|<\varepsilon, |r(x-a)-v|<\varepsilon \},
\end{equation}
such vectors $v$ are called tangent vectors of $S$ at $a$. 

\medskip\noindent

In Section 7, we establish

\medskip\noindent

{\bf Theorem}.\,\, {\it Let $k$ be an integer with $0\leq k\leq r$, let $M_{0}$ be a $n \times n$ skew symmetric matrix of rank $2k$, then the tangent cone of ${\bf C}(n,2r)$ at $M_{0}$ is isometric to ${\bf C}(n-2k,2r-2k) \times  \textbf{R}^{k(2n-2k-1)}$}.

\medskip\noindent

The regular cones ${\bf C}(n,2)$, $(n\geq 4)$, are cones over orientated real Grassmannians, many people has studied their geometries and area-minimizations (see Section 2). In special, the first one ${\bf C}(4,2)$ is just cone over products of two-spheres, it is non-minimizing. By using this fact, we have a corollary for the above theorem.  

\medskip\noindent

{\bf Corollary}.  {\it For the 'nearly regular' points of Pfaffian hypersurfaces ${\bf C}(2n,2n-2)$, their tangent cones are not area-minimizing, where 'nearly regular' means those points \emph{(}matrices\emph{)} have rank $2n-4$.}

\medskip\noindent

 In GMT, an related knowledge states that (cf. Theorem 5.4.3 in \cite{Fe69}) every oriented tangent cone at interior points of an area-minimizing locally rectifiable current $M$ is area-minimizing too. Here, \textit{interior points of $M$} means \textit{points in the support of $M$ but not in the support of $\partial M$}. So, one can conclude that ${\bf C}(2n,2n-2)$ is not area-minimizing if there exists non-minimizing oriented tangent cone at its interior points.

\medskip\noindent

However, it seems that the concrete meanings of tangent cones and \textit{oriented tangent cones} (cf. 9.7 in \cite{M16}) are different. In general, the support of an oriented tangent cone is contained in the support of tangent cone, but equality may not hold. Figure 9.7.2 and Exercise 9.8 in \cite{M16} gives an such example, there each nonnegative axis is an oriented tangent cone, but the tangent cone is the first quadrant.

\medskip\noindent

Then, in the lack of conclusive evidence for the uniqueness of oriented tangent cones for these algebraic varieties, we prefer not to conclude that ${\bf C}(2n,2n-2)$ are definitely non-minimizing after the above corollary.
 
\medskip\noindent

{\bf  Remark 1.2}. The minimal determinantal varieties and Pfaffian varieties are related to the \emph{minimal cubic cone question} posed in \cite{Hs67} (see \cite{PX93}, \cite{T10-1}, \cite{T10-2}, \cite{T19}, \cite{HLT17} for its further researches), there Hsiang find two irreducible, invariant minimal cubic cones: the first $SO(4)$-invariant cubic cone is just the determinantal hypersurface $C(3,3,2)$ (see \cite{T10-2}); the second $SU(4)$-invariant cubic cone is the Pfaffian hypersurface ${\bf C}(6,4)$ as presented in \cite{HLT17}. In addition, there also exists a $Sp(4)$-invariant cubic cone (\cite{Liu89}), and it is interesting to investigate the geometric meanings and area-minimizations of those invariant cones.

\medskip\noindent
\medskip\noindent

\textbf{1.4 Organizations of the paper.} In Section 2, we discuss the definition of Pfaffian varieties, its rank and orbit stratifications. We also discuss the geometries of the regular cones: ${\bf C}(n,2)$, then prove that all Pfaffian varieties are orientable. In Section 3, we briefly explain Lawlor's \textit{directed slicing} method, the \textit{normal space slicings} and the computing formulas for weighting functions. In Section 4, we firstly describe the primary slicings, prove that any two primary slicings must be disjoint or they coincide, compute the primary  weighting functions and the second fundamental forms. In Section 5, we describe the secondary slicings of Kerckhove and Lawlor, then give a low bound estimation for composition weighting functions. In Section 6, we give the proof of the main theorem which is inspired by the original proof given in \cite{KL99}. In Section 7, we study the tangent cones of Pfaffian varieties. 

\medskip\noindent
\medskip\noindent

\section{Pfaffian varieties}

\medskip\noindent



Determinantal varieties $C(p,q,r)$ are the common zero locus of determinants of all $(r+1)\times (r+1)$ submatrices, its skew-symmetric analogs ${\bf C}(n,2r)$ are in fact defined by a family of Pfaffians, to see this, let us recall some works on Pfaffian varieties in algebraic geometry (cf. \cite{He69}, \cite{DP76}, \cite{AD80}, \cite{R00}). 

Let $M$ denote a skew-symmetric matrix of order $n$, its determinant and Pfaffian are denoted by ${\rm det} (M)$ and ${\rm Pf}(M)$, it is well known that ${\rm det}(M)=({\rm Pf}(M))^2$, ${\rm Pf}(M)=0$ when $n$ is odd. Generally, the symmetric algebra $S(\wedge^{2}{\bf R}^{n})$ can be identified with polynomial ring ${\bf R}[x_{ij}]$, $x_{ij}$ are entries of $M$.  Given $2r$ increasing indices $1\leq i_{1}< \ldots < i_{2r}\leq n$, let $[i_{1},\ldots,i_{2r}]$ denote the Pfaffian of the principal minor of $M$ with row and column indices $i_{1},\ldots,i_{2r}$, i.e.
\begin{equation*}
	[i_{1},\ldots,i_{2r}]=\sum_{I_{\sigma}} {\rm sgn}(I_{\sigma})M_{j_{1}k_{1}}\cdots M_{j_{r}k_{r}}
\end{equation*}
where $I_{\sigma}: \{i_{1},\ldots,i_{2r}\}\rightarrow \{j_{1},k_{1},\ldots,j_{r},k_{r}\}$ are permutations satisfying  $j_{t}\leq k_{t} (1\leq t\leq r) $ and $j_{1}\leq \cdots \leq j_{r} $.

De Concini and Procesi (\cite{DP76}) proved that ${\bf R}[x_{ij}]$ is generated by these Pfaffians of principal minors. Then, the common zero locus of determinants of all $(2r+1)$ order and $(2r+2)$ order submatrices of $M$ is reduced to the common zero locus of Pfaffians of principal $(2r+2)$ order minors of $M$ (also see \cite{He69}, \cite{AD80}). Hence, it seems more accurate to call skew-symmetric determinantal varieties--Pfaffian varieties. 
\medskip\noindent
\medskip\noindent

{\bf Definition 2.1}\,\,  For $0\leq r\leq [\frac{n}{2}]-1$, Pfaffian varieties ${\bf C}(n,2r)$ are defined to be the subset of $n$ by $n$ real skew-symmetric matrices of rank at most $2r$, the space of $n$ by $n$ skew-symmetric matrices can be identified with the Euclidean space: $Skew(n,\textbf{R})\cong \textbf{R}^{\frac{n(n-1)}{2}}$, the metric (or the norm) is given by:
\begin{equation}\notag
|M|^2=\frac{1}{2} tr(MM^{T}), M\in Skew(n,\textbf{R}).
\end{equation}
We also call a point(matrix) $N$ of ${\bf C}(n,2r)$ \textit{nearly regular} if the rank of $N$ is equal to $2r-2$.

\medskip\noindent
\medskip\noindent

Similar to determinantal varieties considered in \cite{KL99},  Pfaffian varieties also admit stratifications by the rank:
\begin{equation}\notag
\{0\}={\bf C}(2n,0)\subset {\bf C}(2n,2)\subset \cdots \subset {\bf C}(2n,2n-2),
\end{equation}
the last one: ${\bf C}(2n,2n-2)$ is a singular hypersurface in $Skew(2n,\textbf{R})$ defined by the zero level set of Pfaffian of $2n$ order matrices, we call it \textit{Pfaffian hypersurface},
and 
\begin{equation}\notag
\{0\}={\bf C}(2n+1,0)\subset {\bf C}(2n+1,2)\subset \cdots \subset {\bf C}(2n+1,2n-2),
\end{equation}
the last one has codimension 3 in $Skew(2n+1,\textbf{R})$, and each cone in the above two sequences is the singular set of the next. 

{\bf Remark 2.2}\,\, {\it There are minor differences for the definitions of Pfaffian varieties in algebraic geometry. A. Kuznetsov \cite{Ku06} called Pfaffian varieties the collection of highest rank skew-forms(in dual Euclidean spaces), that only correspond to ${\bf C}(2n,2n-2)$, ${\bf C}(2n+1,2n-2)$ here, and he called other elements in the above two sequences: \textit{generalized Pfaffian varieties}}.

\medskip\noindent

We introduce some notions which will frequently be used throughout this paper. Let $X_{ij}(i<j)$ denote the $n$ by $n$ skew-symmetric matrix with $1$ in the $(i,j)$ position, -$1$ in the $(j,i)$ position and zeros everywhere, denote
$$
M(x_{1},\ldots,x_{k})=\sum_{i=1}^{k}x_{i}X_{2i-1,2i} \in {\bf C}(n,2r),
$$
where $k\leq r$ and $x_{1}\geq \cdots \geq x_{k}$ are nonascending positive number.

Let $E\subset {\bf C}(n,2r)$ be the set of all such $M(x_{1},\ldots,x_{r})$, and set
\begin{equation}\notag
\tilde{E}=\{M(x_{1},\ldots,x_{r}):x_{1}>\cdots > x_{r} > 0 \}.
\end{equation}

\medskip\noindent

The group $SO(n)$ has an adjoint action on the stratified surfaces ${\bf C}(n,2r)$, every orbit intersects $E$ at some $M(x_{1},\ldots,x_{k})(k\leq r)$ with $x_{1}\geq \cdots \geq x_{k}$, and $x_{1},x_{1},\ldots,x_{k},x_{k}$ are the nonascending nonzero singular values. It follows that the closure of $\tilde{E}$: $Clos \ \tilde{E}$ is the orbit space of ${\bf C}(n,2r)$ under adjoint actions, and ${\bf C}(n,2r)$ is the union of orbits. We note here, orbit stratifications are refinements of the rank stratifications, since there exists some orbits through repeated singular values}.

\medskip\noindent

The most singular orbit ${\bf C}(n,2)$ are regular cones over the link:
$$
L=P \cdot  diag \left\lbrace  \begin{pmatrix}
0 & 1\\
-1 & 0
\end{pmatrix},0,\ldots,0 \right\rbrace  \cdot P^{T}, P\in SO(n),
$$
$L$ is an adjoint orbit through base point $X_{12}$. It is isolated by checking the image of orbit spaces: $Clos \ \tilde{E}$(restricted on $x_{1}^2+\cdots+x_{r}^2=1$), hence be minimal in the unit sphere of $Skew(n,\textbf{R})$ by a theorem in \cite{HsL71}, therefore ${\bf C}(n,2)$ is minimal in $Skew(n,\textbf{R})$ as well.

\medskip\noindent

Additionally, the area-minimizations of ${\bf C}(n,2)$ was well studied, let's talk more about them.

\medskip\noindent

The adjoint action of Lie group $SO(n)$ can be realized as the isotropy action of symmetric spaces  $(G,K)=(SO(n)^2,SO(n))$ (often divided into $B$-type and $D$-type in Lie group and symmetric space theories). In this point view, $L$--the adjoint orbit through specially chosen base point, is a symmetric R-space (\cite{Ka02}). And Kanno \cite{Ka02} proved that: except for ${\bf C}(4,2)$, all other cones ${\bf C}(n,2)(n\geq 5)$ are area-minimizing,  based on Lawlor's curvature criterion.

By verifying the Pfaffian of a $4$ by $4$ skew symmetric matrix, after an orthogonal transformation, ${\bf C}(4,2)$ turns out to be the cone over the product of two-spheres:
$$
{\bf C}(4,2)=C(S^2(1)\times S^2(1)),
$$
which cannot be area-minimizing since it is a codimension one singular surface in $\textbf{R}^{6}$ with isolated singularity (even are non-stable, see \cite{Sim68} and \cite{La91}). 
 
\medskip\noindent

The manifold structures of ${\bf C}(n,2)$ are also clear. The real Grassmannians of oriented $m$-planes in $\textbf{R}^{n}$: $\widetilde{G}(m,n;\textbf{R})$ are basic languages in the theory of calibrated geometry (\cite{HL82},\cite{Ha90}), and are naturally embedded into wedge product spaces $\wedge^m \textbf{R}^{n}$ as unit decomposable vectors (cf. \cite{M85}). The  minimality and associated cone's area-minimization were studied in \cite{C88},\cite{JC22}, and the natural embedding of $\widetilde{G}(2,n;\textbf{R})$ into $\wedge^2 \textbf{R}^{n}$ is in fact equivalent to the above symmetric space action (cf. Proposition 5.1 in \cite{JCX22}), i.e. 
$$
{\bf C}(n,2)=the \ cone \ over \ \widetilde{G}(2,n;\textbf{R}).
$$

\medskip\noindent

We conclude this chapter by showing all Pfaffian varieties are orientable, that is different to determinantal varieties (compare with Proposition 4.6 in \cite{KL99}).

\medskip\noindent

{\bf Proposition 2.3}\,\, {\it Pfaffian varieties are all orientable}.

\medskip\noindent

{\bf Proof.} We can only check the orientation of $\Sigma(n,2r)$-the set of skew symmetric matrices of fixed rank $2r$ and all eigenvalues are not equal, it is open and dense in ${\bf C}(n,2r)$. 

$\Sigma (n,2r)$ is the union of adjoint orbits of the same type under the group action of $G=SO(n)$ , its orbit space is  the open, linear cone $\tilde{E}$ in Euclidean space, and the closure of $\tilde{E}$ is the orbit space of ${\bf C}(n,2r)$. Choosing a base point $M_{0}=M(x_{1},\ldots,x_{r})\in \tilde{E}(x_{1}>\cdots>x_{r}>0)$, the isotropy subgroup at $M_{0}$ under adjoint action is checked to be
\begin{align} \notag
H &=S(O(2)\times \cdots \times O(2)\times O(n-2r)) \\ \notag
 &= \left\lbrace {\rm diag} \left\lbrace \varepsilon_{1} R_{\theta_{1}},\cdots,
\varepsilon_{r} 
R_{\theta_{r}},
S \right\rbrace | \varepsilon_{1},\ldots,\varepsilon_{r}\in \pm 1, S\in O(n-2r) \right\rbrace 
\end{align}
where $\prod_{i=1}^r {\rm sgn}(\varepsilon_{i})\cdot {\rm det}S=1 $
and for $i=1,\ldots,r$,
$$
R_{\theta_{i}}=\begin{pmatrix}
\cos \theta_{i} & \sin \theta_{i} \\
-\sin \theta_{i} & \cos \theta_{i}
\end{pmatrix}
$$
represents the $SO(2)$ rotations.

The orientability can be checked by looking at the action of isotropy group $H$ on the tangent space at the identity in $G/H$, and see if it preserves an orientation or not, i.e. checking whether continuous path $PXP^{T}(P\in H,X$ a fixed tangent vectors) preserves the oriented tangent basis of $\Sigma(n,2r)$ at $M_{0}$ (see Proposition 4.6 in \cite{KL99}).

 The oriented tangent basis is given as (see Proposition 4.3 and its proof below): 
 the first group
 $$
\{\ldots, e_{i}, \ldots\}=\{\ldots, X_{2i-1,2i}, \ldots \},1\leq i \leq r,$$ the second group
$$
\{\ldots, f_{ij}, \ldots\}=\{ \ldots,  X_{2i-1,2j-1},X_{2i-1,2j},X_{2i,2j-1},X_{2i,2j}, \ldots\}, 1\leq i<j\leq r,
$$ and the third group 
$$
\{\ldots, g_{i}, \ldots\}=\{\ldots, X_{2i-1,2r+1},\ldots, X_{2i-1,n},X_{2i,2r+1},\ldots,X_{2i,n},\ldots \}, 1\leq i \leq r.$$ 

Denote $R={\rm diag}\{\epsilon_{1} R_{\theta_{1}},\cdots,
\epsilon_{r} 
R_{\theta_{r}}\}$, let $h=\begin{pmatrix}
	R & 0 \\
0 & S
\end{pmatrix}\in H$ be an element of isotropy subgroup, i.e. 
$$
h= \sum_{i=1}^{r}\varepsilon_{i}\left[ cos \theta_{i}(E_{2i-1,2i-1}+E_{2i,2i})+sin \theta_{i}(E_{2i-1,2i}-E_{2i,2i-1})\right]+\sum_{2r+1\leq \alpha,\beta \leq n}s_{\alpha \beta}E_{\alpha\beta},
$$
where $E_{ij}$ represents the matrix with $1$ at $(i,j)$ position, $0$ elsewhere. 

Denote the adjoint actions $h(*)h^{T}$ on the above oriented tangent basis by $\Phi$, concretely, it is given as:
$$
\begin{pmatrix}
R & 0 \\
0 & S
\end{pmatrix}
\begin{pmatrix}
A & B \\
-B^{T} & 0
\end{pmatrix}
\begin{pmatrix}
R^{T} & 0 \\
0 & S^{T}
\end{pmatrix}
=\begin{pmatrix}
RAR^{T} & RBS^{T} \\
-SB^{T}R^{T} & 0
\end{pmatrix}.
$$

Then the following is computed: 
$$\Phi(\ldots, X_{2i-1,2i}, \ldots)=(\ldots, X_{2i-1,2i}, \ldots)\textbf{I}, $$
 it is identity;
$$
\Phi(X_{2i-1,2j-1},X_{2i-1,2j},X_{2i,2j-1},X_{2i,2j})=(X_{2i-1,2j-1},X_{2i-1,2j},X_{2i,2j-1},X_{2i,2j})(\varepsilon_{i}\varepsilon_{j}R_{\theta_{i}}\otimes R_{\theta_{j}}),
$$
where $R_{\theta_{i}}\otimes R_{\theta_{j}}$ denote the Kronecker product of two matrices. Then it is easy to see that the matrices $\varepsilon_{i}\varepsilon_{j}R_{\theta_{i}}\otimes R_{\theta_{j}}$ all have determinants $1$.

And  
$$
\Phi(E_{2i-1,\alpha},E_{2i,\alpha})=(E_{2i-1,\alpha},E_{2i,\alpha})(\varepsilon_{i} R_{\theta_{i}}\otimes S                                                                                                                                                                                                                                                   ),
$$
 where $(E_{2i-1,2r+1},\ldots, E_{2i-1,n},E_{2i,2r+1},\ldots,E_{2i,n})$ is abbreviated to $(E_{2i-1,\alpha},E_{2i,\alpha})$, and ${\rm det}(\varepsilon_{i} R_{\theta_{i}}\otimes S )=(\varepsilon_{i})^{2(n-2r)}\cdot ({\rm det} R_{\theta_{i}})^{n-2r}\cdot ({\rm det}S)^2=1$ is still positive. We note here, by symmetry, that the reader can only check the sign of the action: $RE_{i\alpha}S^{T}$,$1\leq i \leq 2r, 2r+1\leq \alpha \leq n$.

 Therefore, we conclude that ${\bf C}(n,2r)$ are all orientable. $\Box$

\medskip\noindent
\medskip\noindent

\section{Lawlor's directed slicing method}

\medskip\noindent

The directed slicing method was set up in \cite{La98} by Lawlor, then Kerckhove and Lawlor used it to prove determinantal varieties are area-minimizing except two families (\cite{KL99}).

\medskip\noindent

{\bf Definition 3.1 }(cf. Definition 2.1 in \cite{KL99})\,\, Let $M$ be a $k$-dimensional surface in $\textbf{R}^{n}$. For $0\leq d<k$, a \textit{$d$-dimensional slicing} of $M$ is a collection of pairwise disjoint, $d$-dimensional rectifiable subsets of $M$, if the union of these rectifiable subsets covers $M$, $\mathcal{H}^{k}$ a.e., then the slicing is called \textit{full}. 

We are concerned with those $(n-k+d)$-dimensional slicings of $\textbf{R}^{n}$, which its intersections with $M$ consist of $d$-dimensional slicings of $M$, then we call those $d$-dimensional subsets of $M$: \textit{slices}, and those $(n-k+d)$-dimensional subsets of $\textbf{R}^{n}$: \textit{slicing sets} or just \textit{slicings}.

\medskip\noindent

Let $f:\textbf{R}^{n}\rightarrow \textbf{R}^{k-d}$ be a $C^{1}$ map, suppose its restriction on $M$ is also $C^{1}$, a slicing of $M$ given by the level set of $f|_{M}$ is automatically full. The $(k-d)$-dimensional Jacobian of $f$ at $p$ is defined as (see 3.2.1 in \cite{Fe69}, or 3.6 in \cite{M16})
$$
J_{k-d}f={\rm max}||(Df)_{p}(v_{1}\wedge\cdots\wedge v_{k-d})||, 
$$
where $v_{1},\ldots,v_{k-d}$ are unit, orthogonal vectors of $\textbf{R}^{n}$ at $p$, and the $(k-d)$-dimensional Jacobian of $f|_{M}$ at $p\in M$ is:
$$
J_{k-d}(f|_{M})={\rm max}||(Df)_{p}(v_{1}\wedge\cdots\wedge v_{k-d})||, 
$$
where $v_{1},\ldots,v_{k-d}$ are unit, orthogonal tangent vectors of $M$ at $p$.

It can be seen that 
$$
J_{k-d}(f|_{M})\leq (J_{k-d})f|_{M},
$$
since the two Jacobian are defined both on those points $p \in M\subset \textbf{R}^{n}$,  but the tangent vectors for the latter one are chosen from $T_{p}\textbf{R}^{n}$, the equality is attained when every level set of $M$ intersects $M$ orthogonally.

\medskip\noindent

Let $g$ be an integrable function on $M$ with the induced Hausdorff measure, and the area-coarea formula (see 3.2.22 in \cite{Fe69} and Theorem 4.5.1 in \cite{La98}) states that
\begin{equation}\label{areacoarea}
\int_{M}g J_{k-d}(f|_{M})d \mathcal{H}^{k}=\int_{\textbf{R}^{k-d}}\left(\int_{f^{-1}(y)\bigcap M}g d\mathcal{H}^{d}\right) d \mathcal{H}^{k-d}.
\end{equation}

Instead $g$ by the \textit{weighting function} (see 4.6 in \cite{La98}):
$$
w(x)=\frac{1}{(J_{k-d}f)|_{M}}, 
$$
which implies that
$$
{\rm Area} M \geq \int_{M} \frac{1}{(J_{k-d}f)|_{M}} J_{k-d}(f|_{M})d \mathcal{H}^{k}
=\int_{\textbf{R}^{k-d}}\left(\int_{f^{-1}(y)\bigcap M}\frac{1}{(J_{k-d}f)|_{M}} d\mathcal{H}^{d}\right) d \mathcal{H}^{k-d},
$$
the inequality turns out to be equality when every level set $f^{-1}(c)$ intersects $M$ orthogonally.

\medskip\noindent

The inside integral suggests a new weighted area-minimizing question on each slicing set, if we slice $\textbf{R}^{n}$ twice or more, then the inside weighting functions will be the product of the weighting functions produced by iterated slicing every time (see 4.5.3 in \cite{La98}). 

Every composition function $h=g \circ f$ gives the same slicing defined by the level sets, and they are different up a constant on each level set, which will not affect the integral result (see Proposition 4.6.4 and below proof in \cite{La98}).

\medskip\noindent

The slicing of a $C^{1}$ submanifold in $\textbf{R}^{N}$(or $\textbf{R}^{N}$ itself) can be defined by parametrized maps (\cite{La98}). Let $T$ be a properly embedded $n$-dimensional $C^{1}$ submanifold of $\textbf{R}^{N}$, $V\subset \textbf{R}^{n-k+d}$ open, $Z \subset \textbf{R}^{k-d}$ open. 

For a diffeomorphism $F:V\times Z\rightarrow T$, it gives a family of subsets $\{T_{z}\}$ of $T$,
$$
T_{z}:=F(V\times \{z\}), z\in Z,
$$
they are the level sets of $f$ defined by: $f(F(V\times \{z\}))=z,$
thus be a full, $(n-k+d)$-dimensional slicing of $T$. 

The weighting function $w$ for these parametrized slicing sets is computed as follows (see 4.6.2 in \cite{La98}):
\begin{equation}\label{para slicing}
w(F(v,z))=\frac{J_{n}F}{J_{n-k+d}(F|_{V\times \{z\}})}.
\end{equation}

\medskip\noindent

Now, we can introduce a concrete submanifold-like slicing of singular surfaces and its ambient spaces: $normal \ space \ slicings$ (cf. Section 3 in \cite{KL99}).

Given a singular surface $M\subset \textbf{R}^{n}$, near a regular point $p\in M$, $T_{p}\textbf{R}^{n}$ has an orthogonal decomposition: 
$$
T_{p}\textbf{R}^{n}=T_{p}M \oplus N_{p}M, 
$$
where $N_{p}M$ denotes the normal space of $M$ at $p$, it induces a natural extension from slices of $M$ to slicing sets of $\textbf{R}^{n}$ by 'adding' these normal spaces over regular points that belong to some slices of $M$, i.e., 
$$
S_{a}:=\{p+tv|p\in M_{a}, v\in N_{p}M\},
$$
where $p$ is a regular point, $a \in \mathcal{A}$ is an element of some index set, $M_{\mathcal{A}}$ denotes a full slices of $M$, $S_{a}$ denotes the result slicings of $\textbf{R}^{n}$, and we call these $S_{\mathcal{A}}$ normal space slicings. 
    
\medskip\noindent

{\bf Remark 3.2}\,\, {\it 
At any regular point of $M$, those normal space slicing sets should have radii less than the focal distance, which could preserve the positivity of $det(I-t\textbf{A}_{v})$ for any chosen unit normal $v$ in the below; this is a local condition. Globally, every two normal slicing sets should not intersect (see Proposition 4.2), that is similar to the nonintersection of ''normal wedges'' talked in \cite{La91}. For Pfaffian varieties and other singular cones, along radial direction tending to origin, the radii of normal space slicings tend to zero, i.e. they are wedge-shaped (in fact, this can also happen when near other singularities away from origin, check Corollary 4.8 below), and we also call them normal wedges.}

\medskip\noindent

Given a normal space slicing of $M$ that is compatible with local parametrization of $M$ near a regular points $p \in M$, then the weighting functions $w(p+tv)$ and $w(p)$ in \eqref{para slicing} for the normal space slicings satisfy (see  Theorem 3.2 in \cite{KL99}):
\begin{equation}\label{normal function}
w(p+tv)=w(p)\frac{{\rm det}(I-t\textbf{A}_{v})}{J_{d}((I-t\textbf{A}_{v})_{d})},
\end{equation}
where $|t|$ is less than the focal distance of $M$ at a regular point $p$, $p$ belongs to some slice of $M$, and $v$ is a unit normal at $p$. An extra requirement is that $\textbf{A}_{v}$ should be the matrix form of the shape operator w.r.t. an ordered orthonormal basis of $T_{p}M$ such that: \textit{the  first  d  columns  are  tangent  to  the  slices  of  M  through  p}.

\medskip\noindent
\medskip\noindent

\section{Primary slicings and associated weighting functions}

\medskip\noindent

The proof of $\textbf{Main Theorem}$ needs twice slicings. In this chapter, we describe the primary (first) slicings. Keep in mind the following notations introduced in Section 2: $X_{ij},M(x_{1},\ldots,x_{r}),E, \tilde{E}$,
and here we introduce some new ones.

\medskip\noindent

Let $K$ be the subspace spanned by vectors $X_{ij}$ for which $2r<i<j\leq n$, $K$ is the normal space of ${\bf C}(n,2r)$ at regular points.

Let $H=\{M+N:M\in \widetilde{E},N \in K, and \ |N|<M_{2r-1,2r}\}$, where $M_{2r-1,2r}$ denotes the minimum nonzero singular value of $M$, it can be seen as ''adding'' the normal spaces on $\tilde{E}$ bounded by focal radii at every point.

Recall that:
$$\tilde{E}=\{M \in Skew(n,\textbf{R}): M=\sum_{i=1}^{r}x_{i}X_{2i-1,2i},x_{1}>\cdots>x_{r}>0\}.
$$

\medskip\noindent

{\bf Definition 4.1}(compare with Definition 6.1 in \cite{KL99})\,\, The \textit{primary slices} of ${\bf C}(n,2r)$ are defined as the distinct images of adjoint actions on $\tilde{E}$: $P\tilde{E}P^{T}$, $P\in SO(n)$.

The \textit{primary slicing sets} of $Skew(n,\textbf{R})\cong \textbf{R}^{\frac{n(n-1)}{2}}$ are defined as the distinct images of adjoint actions on $H$: $QHQ^{T}$, $Q \in SO(n)$.

The reader can image the shape of primary slicing sets: by Corollary 4.8, they are wedge-shaped near every singularity of ${\bf C}(n,2r)$, including the origin point $O$.
\medskip\noindent

The non-intersection of primary slicing sets is a necessary global condition for applying the method of directed slicings (see Proposition 6.3 in \cite{KL99} and Remark 3.2 here), here we prove that:

\medskip\noindent 

{\bf Proposition 4.2} \,\, {\it Any two primary slicing sets as described above are disjoint or they coincide.}

\medskip\noindent

{\bf Proof.}\,\, Assume there exists two primary slicing set which intersect at $PH_{1}P^{T}=QH_{2}Q^{T}$, $H_{1},H_{2} \in H$, $P,Q\in SO(n)$. By symmetry in position, we can only prove that, for any $H_{3}\in H$, $Q^{T}PH_{3}P^{T}Q\in H$.

Firstly, $PH_{1}P^{T}=QH_{2}Q^{T}$ implies that $H_{1}$ and $H_{2}$ have the same singular values, in special, the same $2r$ largest singular values, so $H_{1}$ and $H_{2}$ have the same $2r$ by $2r$ block in the upper left.

Write $$
H_{1}=\begin{pmatrix}
A & 0 \\
0 & B_{1}
\end{pmatrix}
\ and \
H_{2}=\begin{pmatrix}
A & 0 \\
0 & B_{2}
\end{pmatrix},
$$
where $A=\sum_{i=1}^{r}x_{i}X_{2i-1,2i}(x_{1}>\cdots>x_{r}>0)$, and $B_{1},B_{2}$ are arbitrary skew-symmetric matrices with $|B_{1}|<x_{r},|B_{2}|<x_{r}$.

Moreover, we write 
$$
PH_{1}P^{T}=\left[ P\begin{pmatrix}
L &0 \\ 0 & I
\end{pmatrix}\right]
\begin{pmatrix}
\Sigma & 0\\
0 & B_{1} 
\end{pmatrix}P^{T},
$$
where 
$$
L=diag \left\lbrace \begin{pmatrix}
0 & 1 \\
-1 & 0
\end{pmatrix},\cdots,\begin{pmatrix}
0 & 1 \\
-1 & 0
\end{pmatrix}\right\rbrace 
$$
and 
$$
\Sigma=diag \{ x_{1},x_{1},\cdots,x_{r},x_{r} \}, 
$$ $A=L\Sigma$, similar process for $QH_{2}Q^{T}$.

Then, by the theory of Low Rank Approximation (cf. Proposition 6.2 in \cite{KL99})(we omit the SVD process of $B_{1},B_{2}$ here, which does not affect the result), the unique nearest rank-$2r$ matrix to $PH_{1}P^{T}$ is 
$$
\left[ P\begin{pmatrix}
L &0 \\ 0 & I
\end{pmatrix}\right]
\begin{pmatrix}
\Sigma & 0\\
0 & 0 
\end{pmatrix}P^{T}=P\begin{pmatrix}
A &0 \\ 0 & 0
\end{pmatrix}P^{T},
$$
the uniqueness also implies that
$$
Q^{T}P\begin{pmatrix}
A &0 \\ 0 & 0
\end{pmatrix}P^{T}Q=\begin{pmatrix}
A &0 \\ 0 & 0
\end{pmatrix}.
$$

Denote $Q^{T}P=\begin{pmatrix}
S_{1} & S_{2}\\
S_{3} & S_{4}
\end{pmatrix}\in SO(n)$, then it follows, from above, that $S_{2}=S_{3}=0$ and $S_{1}AS_{1}^{T}=A$.

The $2r$ by $2r$ matrix $A$ only have diagonal blocks of number $r$, view $S_{1}$ as $r$ by $r$ matrix with all entries are $2$ by $2$ blocks, from $S_{1}A=AS_{1}$, notice $x_{1}>\cdots>x_{r}>0$, it concludes that $S_{1}$ only have diagonal blocks too.

Now each block of $S_{1}$ belongs to $O(2)$, then, $$
Q^{T}P=\begin{pmatrix}
S_{1} & 0\\
0 & S_{4}
\end{pmatrix}\in S(O(2)\times \cdots \times O(2)\times O(n-2r))
$$, and finally, $Q^{T}PH_{3}P^{T}Q$ will send every $H_{3}\in H$ back to $H$.  $\Box$ 

\medskip\noindent
\medskip\noindent

In the next, we compute the primary (first) weighting functions, since adjoint actions are isometries, then we can only perform the computations on those points belonging to $\tilde{E}$.

\medskip\noindent

{\bf Proposition 4.3}\,\, {\it Let  $M(\overrightarrow{x})=M(x_{1},\ldots,x_{r})=\sum_{i=1}^{r}x_{i}X_{2i-1,2i} \in \widetilde{E}$ be a regular point of ${\bf C}(n,2r)$, then the first weighting function in \eqref{normal function} of ${\bf C}(n,2r)$ at $M(\overrightarrow{x})$ is given by 
\begin{equation}\notag
w_{1}(M(\overrightarrow{x}))=\prod_{1\leq i<j\leq r}(x_{i}^2-x_{j}^2)^2\prod_{i=1}^{r}x_{i}^{2(n-2r)}.
\end{equation}}

\medskip\noindent


{\bf Proof.} \,\, There is an natural local parametrization for ${\bf C}(n,2r)$ coming from Lie algebra, when seen it as the collection of adjoint orbits: 
\medskip\noindent

Set $M(\overrightarrow{x})=M(x_{1},\ldots,x_{r})=\sum_{i=1}^{r}x_{i}X_{2i-1,2i}$, fix $\overrightarrow{x}$,  let $\overrightarrow{t}=(t_{1},\ldots,t_{r})$, then $M(\overrightarrow{x}+\overrightarrow{t})$ gives an local parametrization of $\tilde{E}$.

Choose the indices $\overrightarrow{s}=(s_{ij}|1\leq i<i+1<j\leq n, i\leq 2r)$, let $P(\overrightarrow{s})={\rm Exp}(\sum_{i,j} s_{ij}X_{ij})$, then the parametrization at $M(\overrightarrow{x})$ is:
$$
F(\overrightarrow{s},\overrightarrow{t})=P(\overrightarrow{s})M(\overrightarrow{x}+\overrightarrow{t})P(\overrightarrow{s})^{T}.
$$

The tangent map at $\overrightarrow{s}=\overrightarrow{0},\overrightarrow{t}=\overrightarrow{0}$ is given by: $F_{*}(\frac{\partial}{\partial s_{ij}})=[X_{ij},M(\overrightarrow{x})] $, $F_{*}(\frac{\partial}{\partial t_{i}})=X_{2i-1,2i}$.

\medskip\noindent








By directly computing, the matrix $DF$ has the following non-zero entries:

\medskip\noindent

For those $j\leq 2r$,

(i): $i=2p<j=2q$,
$$
\frac{\partial F_{i-1,j}}{\partial s_{ij}}=-x_{i/2}, \frac{\partial F_{i,j-1}}{\partial  s_{ij}}=-x_{j/2};
$$

(ii): $i=2p<j=2q-1$,
$$
\frac{\partial F_{i-1,j}}{\partial s_{ij}}=-x_{i/2}, \frac{\partial F_{i,j+1}}{\partial  s_{ij}}=x_{(j+1)/2};
$$

(iii): $i=2p-1<j=2q(p<q)$,
$$
\frac{\partial F_{i,j-1}}{\partial s_{ij}}=-x_{j/2}, \frac{\partial F_{i+1,j}}{\partial  s_{ij}}=x_{(i+1)/2};
$$

(iv): $i=2p<j=2q-1$,
$$
\frac{\partial F_{i,j+1}}{\partial s_{ij}}=x_{(j+1)/2}, \frac{\partial F_{i+1,j}}{\partial  s_{ij}}=x_{(i+1)/2}.
$$

\medskip\noindent

For those $j>2r$, when $i$ is even, $\frac{\partial F_{i-1,j}}{\partial s_{ij}}=-x_{\frac{i}{2}}$; when $i$ is odd, $\frac{\partial F_{i+1,j}}{\partial s_{ij}}=x_{\frac{i+1}{2}}$.


\medskip\noindent

And, $\frac{\partial F_{2i-1,2i}}{\partial t_{i}}=1$.

\medskip\noindent

We divided them into the following three families of block matrices:

$$
\begin{pmatrix}
\frac{\partial F_{12}}{\partial t_{1}} & \frac{\partial F_{34}}{\partial t_{1}} & \cdots & \frac{\partial F_{2r-1,2r}}{\partial t_{1}}\\
\frac{\partial F_{12}}{\partial t_{2}} & \frac{\partial F_{34}}{\partial t_{2}} & \cdots & \frac{\partial F_{2r-1,2r}}{\partial t_{2}} \\
\vdots & \vdots & \vdots & \vdots \\
\frac{\partial F_{12}}{\partial t_{r}} & \frac{\partial F_{34}}{\partial t_{r}} & \cdots & \frac{\partial F_{2r-1,2r}}{\partial t_{r}}
\end{pmatrix}=\textbf{I}.
$$

And the matrices $$
\begin{pmatrix}
\frac{\partial F_{2p-1,2q-1}}{\partial s_{2p-1,2q-1}} &\frac{\partial F_{2p-1,2q}}{\partial s_{2p-1,2q-1}} &\frac{\partial F_{2p,2q-1}}{\partial s_{2p-1,2q-1}} &\frac{\partial F_{2p,2q}}{\partial s_{2p-1,2q-1}}\\
\frac{\partial F_{2p-1,2q-1}}{\partial s_{2p-1,2q}} &\frac{\partial F_{2p-1,2q}}{\partial s_{2p-1,2q}} &\frac{\partial F_{2p,2q-1}}{\partial s_{2p-1,2q}} &\frac{\partial F_{2p,2q}}{\partial s_{2p-1,2q}}\\
\frac{\partial F_{2p-1,2q-1}}{\partial s_{2p,2q-1}} &\frac{\partial F_{2p-1,2q}}{\partial s_{2p,2q-1}} &\frac{\partial F_{2p,2q-1}}{\partial s_{2p,2q-1}} &\frac{\partial F_{2p,2q}}{\partial s_{2p,2q-1}}\\
\frac{\partial F_{2p-1,2q-1}}{\partial s_{2p,2q}} &\frac{\partial F_{2p-1,2q}}{\partial s_{2p,2q}} &\frac{\partial F_{2p,2q-1}}{\partial s_{2p,2q}} &\frac{\partial F_{2p,2q}}{\partial s_{2p,2q}}
\end{pmatrix}
$$
are equal to 
$$
\begin{pmatrix}
0&x_{q}&x_{p}&0\\
-x_{q}&0&0&x_{p}\\
-x_{p}&0&0&x_{q}\\
0&-x_{p}&-x_{q}&0\\
\end{pmatrix}:=G_{p q},
$$
where $1\leq p<q \leq r$.

The third families of matrices
$$
\begin{pmatrix}
\frac{\partial F_{2i-1,2r+1}}{\partial s_{2i-1,2r+1}}& \cdots &\frac{\partial F_{2i-1,n}}{\partial s_{2i-1,2r+1}} &\frac{\partial F_{2i,2r+1}}{\partial s_{2i-1,2r+1}} &\cdots
&\frac{\partial F_{2i,n}}{\partial s_{2i-1,2r+1}}\\
\vdots &\vdots &\vdots &\vdots &\vdots &\vdots \\
\frac{\partial F_{2i-1,2r+1}}{\partial s_{2i-1,n}}& \cdots &\frac{\partial F_{2i-1,n}}{\partial s_{2i-1,n}} &\frac{\partial F_{2i,2r+1}}{\partial s_{2i-1,n}} &\cdots
&\frac{\partial F_{2i,n}}{\partial s_{2i-1,n}}\\
\frac{\partial F_{2i-1,2r+1}}{\partial s_{2i,2r+1}}& \cdots &\frac{\partial F_{2i-1,n}}{\partial s_{2i,2r+1}} &\frac{\partial F_{2i,2r+1}}{\partial s_{2i,2r+1}} &\cdots
&\frac{\partial F_{2i,n}}{\partial s_{2i,2r+1}}\\
\vdots &\vdots &\vdots &\vdots &\vdots &\vdots \\
\frac{\partial F_{2i-1,2r+1}}{\partial s_{2i,n}}& \cdots &\frac{\partial F_{2i-1,n}}{\partial s_{2i,n}} &\frac{\partial F_{2i,2r+1}}{\partial s_{2i,n}} &\cdots
&\frac{\partial F_{2i,n}}{\partial s_{2i,n}}
\end{pmatrix}
$$
are equal to
$$
\begin{pmatrix}
0& x_{i}\textbf{I}\\
-x_{i}\textbf{I}& 0
\end{pmatrix}:=H_{i},
$$
where $1\leq i\leq r$.

\medskip\noindent

After removing invalid variables, it follows from above ordering:
$$
DF= {\rm diag} \left\lbrace \textbf{I},\cdots,G_{pq}(p<q),\cdots,H_{1},\cdots,H_{r} \right\rbrace .
$$

\medskip\noindent

The determinant of $DF$ is the product of these blocks. The $r$-dimensional Jacobian of $F$ which is restricted to slices of ${\bf C}(n,2r)$, $ J_{r}F|_{F^{-1}(slice)} $, is equal to the Jacobian of the formal $r$ columns, which is just $1$. Then, follow \eqref{para slicing}, the primary weighting function on ${\bf C}(n,2r)$ is:
\begin{equation}\label{weighted slices}
\begin{aligned}
w_{1}&={\rm det} \ DF= \prod_{1\leq i<j\leq r} {\rm det} \ G_{ij}\times \prod_{i=1}^{r}{\rm det} \ H_{i}\\
&=\prod_{1\leq i<j\leq r}(x_{i}^2-x_{j}^2)^2\prod_{i=1}^{r}x_{i}^{2(n-2r)}.
\end{aligned}
\end{equation} $\Box$


\medskip\noindent

Recall the formula \eqref{normal function}, to obtain the value of $w_{1}$ on primary slicing sets--those normal slicing sets, we also need to compute the second fundamental forms of ${\bf C}(n,2r)$ in $Skew(n,\textbf{R})\simeq\textbf{R}^{\frac{n(n-1)}{2}}$, and the minimality of ${\bf C}(n,2r)$ is a natural corollary that reprove the result in \cite{Ko21}. Finally, we shall give a lower bound estimation for the first weighting function on primary slicing sets in \eqref{normal function}, and that is enough for further using.

\medskip\noindent

By the isometric action, we could also perform the computations of second fundamental forms at those regular points $M_{0}=\sum_{i=1}^{r}x_{i}X_{2i-1,2i}\in \tilde{E}$. Recall that a basis of tangent vectors of ${\bf C}(n,2r)$ at $M_{0}$ is spanned by $X_{ij}(1\leq i<j\leq n,i\leq 2r)$, and the normal space to the cone at $M_{0}$ is spanned by $X_{ij}(2r<i<j)$, those bases are unit orthogonal vectors under the given metric. The second fundamental forms can be calculated by examining every 2-dimensional section of the cone (cf. p.44 in \cite{Simon83}, and p.415 in \cite{KL99}).

\medskip\noindent

Let$X_{a,h}$ and $X_{b,l}$ be two tangent vectors at $M_{0}$, without loss of generality, assume $a\leq b$, then a key observation is that:

\medskip\noindent

{\bf Lemma 4.4}\,\, {\it Adding a linear combination of $X_{a,h}$ and $X_{b,l}$ to $M_{0}$ raises the rank if and only if $a=2i-1,b=2i$ for some $i=1,\ldots,r$ and $h>2r,l>2r$}.

{\it Proof.}\,\, It is a direct consequence when discuss separately all the following cases: $a=2i<b=2j$,$a=2i<b=2j-1$,$a=2i-1<b=2j(i<j)$ and $a=2i-1,b=2i$. $\Box$
\medskip\noindent

Therefore, for any pair $X_{2i-1,h}$ and $X_{2i,l}$, there exists a curve $f(y_{1},y_{2})\in {\bf C}(n,2r)$ in a small neighborhood of $M_{0}$:
$$
f(y_{1},y_{2})=M_{0}+y_{1}X_{2i-1,h}+y_{2}X_{2i,l}+\frac{y_{1}y_{2}}{x_{i}}X_{h,l},
$$
it satisfies $\frac{\partial f}{\partial y_{1}}(0,0)=X_{2i-1,h}$, $\frac{\partial f}{\partial y_{2}}(0,0)=X_{2i,l}$, the associated second fundamental forms are:
$$
\overrightarrow{\textbf{B}}(X_{2i-1,h},X_{2i,l})=-\left[ \frac{\partial^{2}f}{\partial y_{1} \partial y_{2}}(0,0)\right]^{\perp} =-\frac{1}{x_{i}}X_{h,l},
$$
where $i=1,\ldots,r$ and $h,l>2r$, since $X_{h,l}$ is normal.
\medskip\noindent

For any other pairs $X_{a,h}$ and $X_{b,l}$, those curves (small pieces) can be just chosen: $M_{0}+y_{1}X_{a,h}+y_{2}X_{b,l}$, which implies that the associated second fundamental forms are zeros. 

\medskip\noindent

We give an order of tangent vectors $X_{ij}(1\leq i<j\leq n,i\leq 2r)$ as follows:

$$
X_{12},\ldots,X_{2r-1,2r},X_{13},\ldots,X_{1,2r},X_{23},\ldots,X_{2,2r},\ldots,X_{2r-2,2r};
$$
$$
X_{1,2r+1},\ldots,X_{1,n},X_{2,2r+1},\ldots,X_{2,n},\ldots,X_{2r,2r+1},\ldots,X_{2r,n}.
$$

The formal $r$ vectors are tangent to the primary slices $\tilde{E}$ of ${\bf C}(n,2r)$. Under this basis, the second fundamental form at points of $\tilde{E}$ are given by:

\medskip\noindent

{\bf Proposition 4.5}\,\, {\it Let $M_{0}=\sum_{i=1}^{r}x_{i}X_{2i-1,2i}\in \tilde{E}$, let $v=\begin{pmatrix}0 & 0\\ 0 & B \end{pmatrix} $ be a unit normal vector to ${\bf C}(n,2r)$ at $M_{0}$, i.e. $|B|=\sqrt{ \frac{1}{2} {\rm tr}BB^{T}}=1$. With the above ordered basis, the associated matrix of shape operator $\textbf{A}_{v}$ has the block forms $\textbf{A}_{v}=\begin{pmatrix}
0 & 0\\ 0 & C \end{pmatrix},$
where 
\begin{equation}\label{2nd}
C=diag\{L_{1},\ldots,L_{r}\}, L_{i}=\begin{pmatrix}
0 & \frac{1}{x_{i}}B^{T}\\
\frac{1}{x_{i}}B & 0
\end{pmatrix},
\end{equation}
for $1\leq i \leq r$.
}

\medskip\noindent

The shape operator $\textbf{A}_{v}$ has trace zero, hence

{\bf Corollary 4.6} (The minimality of Pfaffian varieties, also see \cite{Ko21})\,\, {\it The Pfaffian varieties ${\bf C}(n,2r)$ are minimal singular surfaces in  $Skew(n,\textbf{R})$ under the natural inclusion}. 

\medskip\noindent

The following lemma is crucial for estimating the weighting functions on primary slicing sets, the proof is omitted.

{\bf Lemma 4.7} (A square version of Lemma 6.8 in \cite{KL99}) \,\, {\it 
Let $S$ be a positive, semidefinite, symmetric matrix with trace $2$ and repeated nonnegative eigenvalues $\lambda_{1},\lambda_{1},\ldots,\lambda_{n},\lambda_{n}$. Let $\lambda$ be the largest eigenvalues of $S$, and let $\tau$ be a real number with $0\leq \tau \leq \frac{1}{\lambda}$. Then, $det(I-\tau S)\geq (1-\tau)^2$, the equality is attained if and only if $t=0$ or $\lambda=1$, i.e., regardless of the ordering, the eigenvalues are $1,1,0,0,\cdots,0,0$.}

\medskip\noindent 

Now, \eqref{2nd} implies that
$$
{\rm det}(I-t\textbf{A}_{v})=\prod_{i=1}^{r}{\rm det}\left(I-\frac{t^2}{x_{i}^2}BB^{T}\right), 
$$
note that $BB^{T}$ has pairs of repeated eigenvalues. By the above lemma, it deduces that

\begin{equation}\notag
{\rm det}\left(I-\frac{t^2}{x_{i}^2}BB^{T}\right)\geq \left(1-\frac{t^2}{x_{i}^2}\right)^{2},
\end{equation}
the equality is attained if and only if $t=0$, or the eigenvalues of $BB^{T}$ are $1,1,0,0,\cdots,0,0$. The later condition is equivalent to ${\rm rank}\ B=2$, which also means that $B\in \widetilde{G}(2,n-2r;\textbf{R})$ in the sense of Section 2 since $v$ is a unit vector.

\medskip\noindent

Then we have a corollary for focal radius:

{\bf Corollary 4.8} \,\, {\it The focal radius of ${\bf C}(n,2r)$ at $M_{0}=M(x_{1},\ldots,x_{r})$ $\in \tilde{E}$ is the smallest singular values $x_{r}$, which tend to zero when near singularities. Therefore, the slicing sets are wedge-shaped.} 

\medskip\noindent

Combined with \eqref{weighted slices}, we obtain the main result of this chapter:

{\bf Proposition 4.9} \,\, {\it Let $M_{0}=\sum_{i=1}^{r}x_{i}X_{2i-1,2i}\in \tilde{E}$, let $v$ be unit normal to ${\bf C}(n,2r)$ at $M_{0}$, and let $t$ be a real number with $|t|<x_{r}$. Then the weighting function $w_{1}$ for the primary slicing of $\textbf{R}^{\frac{n(n-1)}{2}}$ satisfies
\begin{equation}\label{primary weighted function normal wedge}
w_{1}(M(\overrightarrow{x})+tv)\geq \prod_{1\leq i <j \leq r}(x_{i}^2-x_{j}^2)^2\prod_{i=1}^{r}x_{i}^{2(n-2-2r)}\prod_{i=1}^{r}(x_{i}^2-t^2)^2,
\end{equation}
the equality is attained if and only if $t=0$ or
$$
v=\begin{pmatrix}0 & 0\\0 & B\end{pmatrix},B\in {\bf C}(n-2r,2).$$with ${\rm tr}BB^{T}=2$. In special, the equality is naturally satisfied for Pfaffian hypersurfaces ${\bf C}(2n,2n-2)$.}

\medskip\noindent
\medskip\noindent

\section{The secondary slicings and composite weighting functions}

\medskip\noindent

The factors $(x_{i}^2-t^2)$ in \eqref{primary weighted function normal wedge} give a motivation to define the secondary slicings which are 1-dimensional, they will decompose the wedge-shaped primary slicing sets into disjoint hyperbolic curves (cf. Figure 2 in \cite{La98}).  

\medskip\noindent

{\bf Definition 5.1} (see Chapter 7 in \cite{KL99}) \,\, Let $H,K$ be as before, let $M(\overrightarrow{x})=M(x_{1},\ldots,x_{r})=\sum_{i=1}^{r}x_{i}X_{2i-1,2i}\in \tilde{E}$, let $v\in K$ be a unit normal at $M(\overrightarrow{x})$. Then for $|t|<x_{r}$,  define $h=(h_{1},\ldots,h_{r}):H\rightarrow \textbf{R}^{r}$ as:
$$
h_{i}(M(\overrightarrow{x})+tv)=\frac{x_{i}^2-t^2}{2},i=1,\ldots,r,
$$
and the secondary slicings of $H$ are defined to be distinct images of the level sets of $h$, denoted them by $H_{c}$,
$$
H_{c}:=h^{-1}(c),
$$
where $c=(c_{1}^2/2,\ldots,c_{r}^2/2)$. Moreover, the union of secondary slicings sets $H_{c}$ covers $H$ a.e., and  every $H_{c}$ intersect orthogonally with $\tilde{E}$ at a single point.

\medskip\noindent

The weighting function for the secondary slicings is computed in \cite{KL99} to be

\begin{equation}\label{second weighting function}
w_{2}(M(\overrightarrow{x})+tv)=\frac{1}{\prod_{i=1}^{r}x_{i}\sqrt{1+t^2(\sum_{i=1}^{r}\frac{1}{x_{i}^{2}}})}
\end{equation}

\medskip\noindent
\medskip\noindent

Now, we can show that 

{\bf Proposition 5.2} \,\, {\it When $n-2r\geq 3$, the minimum value of the composite weighting function $w_{1}w_{2}$ within each secondary slicing set equals to its value taken at ${\bf C}(n,2r)$}.

{\bf Proof.} \,\, Combined \eqref{primary weighted function normal wedge} and \eqref{second weighting function}, $w_{1}w_{2}$ is bounded below by 
\begin{equation}
w_{1}w_{2}(M(\overrightarrow{x})+tv)\geq \frac{\prod_{1\leq i <j \leq r}(x_{i}^2-x_{j}^2)^2\prod_{i=1}^{r}x_{i}^{2n-4r-5}\prod_{i=1}^{r}(x_{i}^2-t^2)^2}{\sqrt{1+t^2(\sum_{i=1}^{r}\frac{1}{x_{i}^{2}})}},
\end{equation}
the equality is attained if and only if $t=0$ or $v$ is a unit normal vector in matrix form $$v=\begin{pmatrix}0 & 0\\0 & B\end{pmatrix},B\in {\bf C}(n-2r,2).$$

\medskip\noindent

Along every fixed secondary slicing set $H_{c}(c=(c_{1}^2/2,\ldots,c_{r}^2/2))$, the factors $(x_{i}^2-x_{j}^2)$ and $(x_{i}^2-t^2)$ remain unchanged, then the proof is done if (not only if!) we can show that
$$
\frac{\prod_{i=1}^{r}(c_{i}^2+t^2)^{\frac{2n-4r-5}{2}}}{\sqrt{1+t^2(\sum_{i=1}^{r}\frac{1}{c_{i}^{2}+t^2}})}\geq \prod_{i=1}^{r}(c_{i}^2)^{\frac{2n-4r-5}{2}},
$$
and it is equivalent to  
$$
\prod_{i=1}^r\left(c_i^2+t^2\right)^{2n-4r-4} \geqslant \prod_{i=1}^r\left(c_i^2\right)^{2n-4r-5}\left\{\sum_{j=1}^r t^2 \prod_{k \neq j}\left(c_k^2+t^2\right)+\prod_{j=1}^r\left(c_j^2+t^2\right)\right\}.
$$

\medskip\noindent

We apply lemma 8.2 in \cite{KL99}, if and only if the integer $2n-4r-4\geq 2$, i.e.  $n-2r\geq 3$, the above formula is satisfied. $\Box$.

Then, the directed slicing method cannot be applied to the case of hypersurfaces: ${\bf C}(2n,2n-2)$, we will talk about their tangent cone behaviors in Section 7.
\medskip\noindent
\medskip\noindent

\section{Proof of \textbf{Main Theorem}}

\medskip\noindent

\textbf{Main Theorem}. {\it Pfaffian varieties are stratified area-minimizing cones except for the case of hypersurfaces}.

\medskip\noindent

Inspired by the original proof given in \cite{KL99}, here we give the proof of \textbf{Main Theorem}. To make it clear, we isolate the following proposition.  

\medskip\noindent

Let $C_{0}$ be the truncated part of ${\bf C}(n,2r)$ inside the unit ball of $Skew(n,\textbf{R})\cong\textbf{R}^{N},N=\frac{n(n-1)}{2}$, let $S$ be a rectifiable current with the same boundary as $C_{0}$. Those distinct rotated images of $H$: $QHQ^{T},Q\in SO(n)$, the primary slicing sets of $\textbf{R}^{N}$, are denoted by $T_{a}(a\in \mathcal{A})$, and the associated primary slices of $C_{0}$ and $S$ are denoted by $C_{a}$ and $S_{a}$: $C_{a}=T_{a}\cap C_{0}$, $S_{a}=T_{a}\cap S$.  
\medskip\noindent

{\bf Proposition 6.1} (see Theorem 8.3 in \cite{KL99}) \,\, {\it 
If $n-2r\geq 3$, then every primary slice $C_{a}$ solves a special \textit{weighted Plateau's problem} where the comparison surfaces are taken to be the intersection of $T_{a}$ with an arbitrary rectifiable current $S$, and $\partial S= \partial C_{0}$, the weighting function is taken to be the first weighting function in \eqref{normal function}, i.e. 
$$
\int_{C_{a}}w_{1}\leq \int_{S_{a}}w_{1},
$$
where $C_{a}=T_{a}\cap C_{0}$, $S_{a}=T_{a}\cap S$.
 }

\medskip\noindent

{\bf Proof.} \,\, 
Since the adjoint actions of $SO(n)$ on ${\bf C}(n,2r)$ are isometries and $S$ is chosen arbitrarily, we could only prove that
$$
\int_{C_{0}\cap H}w_{1}\leq \int_{S\cap H }w_{1}.
$$

\medskip\noindent

Let $U$ be the union of the secondary slicing sets $H_{c}$ such that $H_{c}$ intersects ${\bf C}(n,2r)$ within the open unit ball. Clearly, $U$ is open in $H$ and covers $C_{0}\cap H$ almost everywhere, by the area-coarea formula, note that every slicing set $H_{c}\subset U$ intersects $C_{0}$ orthogonally, we have

$$
\int_{C_{0}\cap H}w_{1}=\int_{C_{0}\cap U}w_{1}= \int_{c \in \textbf{R}^{r}}w_{1}w_{2}(H_{c}\cap C_{0}) 
$$

The slicing commutes with the boundary operator, up to sign(see 4.3 in \cite{Fe69} and 4.11 in \cite{M16}). Now $S-C_{0}$ bounds an higher dimensional current, denoted by $R$, i.e., $\partial R=S-C_{0}$, then 
$$
(S-C_{0})\cap H_{c}= (\partial R) \cap H_{c} = \pm \partial (R \cap H_{c}),
$$
by the dimension argument, $R \cap H_{c}$ is a curve or set of curves, and its boundary will consists at least two points, which implies that $S$ intersects almost every slicing set $H_{c}\subset U$. Now, combined with Proposition 5.2, it deduces that

$$
\int_{c \in \textbf{R}^{r}}w_{1}w_{2}(H_{c}\cap C_{0})\leq 
\int_{c \in \textbf{R}^{r}} \sum_{x\in H_{c}\cap S} w_{1}w_{2}(x) \leq \int_{S\cap H}w_{1},
$$
the first inequality follows form Proposition 5.2, the last inequality is again by the area-coarea formula, the sum indices arise from the fact that every $H_{c}$ intersects $S$ at least one point, maybe more! 

Then 
$$
\int_{C_{0}\cap H}w_{1} \leq \int_{S\cap H}w_{1}.
$$
$ \Box $

\medskip\noindent

\begin{figure}[hbt]
\begin{center}
  \includegraphics{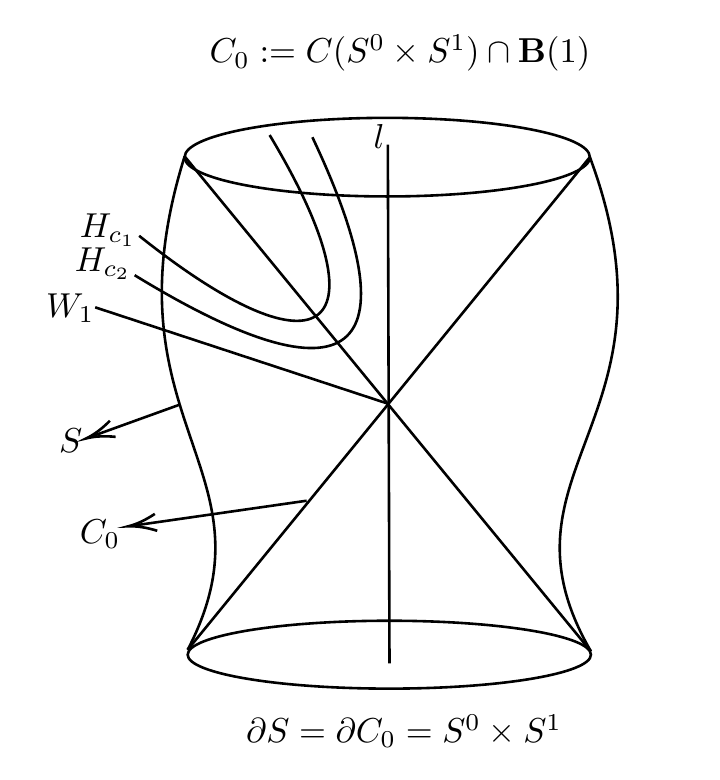}
  \caption{Directed slicings for $C(S^0 \times S^1)$ }
  \end{center}
\end{figure}

Figure 1 gives a sketchy pictured description for Lawlor's directed slicing method on $C(S^0\times S^1) \subset \textbf{R}^3$(obviously, it is not area-minimizing). In figure 1, $l$ together with $W_{1}$ bound a normal wedge of $C_{0}$(there exists another normal wedge in the low part of $C_{0}$), and this normal wedge is further sliced up by parabolic secondary slicing sets $H_{c_{1}},H_{c_{2}},\ldots$, we can see that the comparison surface $S$ intersects every $H_{c}$ at least one point, but part of $S$ are not covered by $T$.

\medskip\noindent

{\bf Proof of the Main Theorem} \,\, Let $S$ be a rectifiable current with the same boundary as $C_{0}$, follow Proposition 6.1, note that $T$ covers $C_{0}$ a.e., then  
$$
{\rm Area}(C_{0})={\rm Area}(C_{0}\cap T)=\int_{a\in \mathcal{A}}\int_{C_{a}}w_{1}\leq \int_{a\in \mathcal{A}}\int_{S_{a}}w_{1} \leq {\rm Area}(S \cap T) \leq {\rm Area}(S),
$$ 
the second equality and the second inequality are both by the area-coarea formula. The final inequality follows the fact that generally $T$ do not cover comparison surface $S$ a.e., see above figure. This completes the proof of the \textbf{Main Theorem}. $\Box$

\medskip\noindent
\medskip\noindent

\section{Tangent cones of Pfaffian varieties}

\medskip\noindent





{\bf Definition 7.1} (Tangent cone, see 3.1.21 in \cite{Fe69}).\,\,  For an normed vector space $X$, a subset $S\subset X$ and a point $a\in X$, we define the tangent cone of $S$ at $a$:
\begin{equation}\notag
{\rm Tan}(S,a)=\{v\in X: \forall \varepsilon >0, \exists x \in S  \ \& \ r\in \textbf{R}^{+}, s.t. \ |x-a|<\varepsilon, |r(x-a)-v|<\varepsilon \},
\end{equation}
such vectors $v$ are called tangent vectors of $S$ at $a$. 

\medskip\noindent

The computation of tangent cones of ${\bf C}(n,2r$) mainly depends on the distance behavior of skew-symmetric matrices of different ranks, which is often known as the theory of Low Rank Approximation in matrix analysis. By studying these distance functions,  we establish a similar result to determinantal varieties (compare with Proposition 4.2 in \cite{KL99}).

\medskip\noindent

{\bf Theorem 7.2} (Tangent cone to Pfaffian varieties).\,\, {\it Let $k$ be an integer with $0\leq k\leq r$, let $M_{0}$ be a $n$ by $n$ skew symmetric matrix of rank $2k$, without losing generality by the isometric adjoint action, we can set $M_{0}=M(x_{1},\ldots,x_{k},0,\ldots,0)=\sum_{i=1}^{k}x_{i}X_{2i-1,2i}$.
Let $V$ be a $n$ by $n$ skew symmetric matrix having block form
$$
V=\begin{pmatrix}
A & B\\
-B^{T} & D
\end{pmatrix},
$$
where the block $A $ is $2k$ by $2k$, then 

\begin{equation}\label{tangent vector}
V\in {\rm Tan}({\bf C}(n,2r),M_{0}) \Leftrightarrow {\rm rank} \ D\leq 2r-2k.
\end{equation}

This implies that
\begin{equation}\label{tan cone}
{\rm Tan}({\bf C}(n,2r),M_{0})\cong {\bf C}(n-2k,2r-2k) \times \textbf{R}^{k(2n-2k-1)}.
\end{equation}}

{\bf Remark 7.3} {\it Following \eqref{tan cone}, the tangent cones of ${\bf C}(n,2r)$ at $M_{0}$ are cylindrical, if $M_{0}$ is a nearly regular point for the meaning that its rank is $2r-2$, then the cross section is a regular cone:
$$
{\rm Tan}({\bf C}(n,2r),M_{0})\cong {\bf C}(n-2r+2,2) \times \textbf{R}^{(r-1)(2n-2r+1)}\cong C(\widetilde{G}(2,n-2r+2;\textbf{R})) \times \textbf{R}^{(r-1)(2n-2r+1)},
$$
thus gives a realization of tangent cones for the oriented real Grassmannian manifolds.

\medskip\noindent

The condition $2r-2k=2,n-2k\geq 5$ in \eqref{tan cone} and the condition $n-2r+2\geq 5$ both are equivalent to $n-2r\geq 3$ except trivial cases, i.e. when $n-2r\geq 3$, the above tangent cones are area-minimizing,  that coincide with the {\bf Main Theorem}.}

\medskip\noindent

{\bf Proof of Theorem 7.2.}\,\, We are aim to prove \eqref{tangent vector}, the proof is an interpretation of the original one given in \cite{KL99} for determinantal varieties, and some differences are also exhibited.

First, assume ${\rm rank} \ D\leq 2r-2k$, we will show that for any chosen $\varepsilon >0$, there exists $X\in {\bf C}(n,2r)$, s.t. 
\begin{equation}\label{assume one}
|X-M_{0}-tV|=O(t^2),
\end{equation}
it follows, for some $t_{0}$, when $0<t<t_{0}$, 
$$
|X-M_{0}|\leq t|V|+O(t^2)<\varepsilon 
$$
and
$$
\left|\frac{1}{t}(X-M_{0})-V\right|=O(t)<\varepsilon,
$$
when $0<t<t_{1}$ for some  $t_{1}$.

If we choose $r=\frac{1}{t}$ as in the definition of the tangent cone, where $t\in (0, min\{t_{0},t_{1}\})$. Then, it is easy to see that $V$ is a tangent vector at $M_{0}$.

Denote the upper left of $M_{0}$ by $M$, and choose $X$ as
$$
X=\begin{pmatrix}
M+tA & tB\\
-t B^{T} & tD-t^2 B^{T}(M+tA)^{-1}B
\end{pmatrix},
$$
then by block Gaussian elimination, we can see that \eqref{assume one} is satisfied.

Conversely, assuming ${\rm rank} \ D>2r-2k$, we will show that for any $t$, there exists some $X_{t}\notin {\bf C}(n,2r)$ with 
\begin{equation}\label{assume two}
\begin{cases}
|X_{t}-(M_{0}+tV)|& = O(t^2),\\
{\rm Dist}(X_{t},{\bf C}(n,2r))&\geq c \ t,
\end{cases}
\end{equation}
for some constant $c$.

Then, it follows that for any $ Y\in {\bf C}(n,2r)$, $|Y-(M_{0}+ t V)| \geq \widetilde{c} \ t$, where $\widetilde{c}$ is another constant, i.e. 
$$
\left|\frac{1}{t}(Y-M_{0})-V\right|\geq \widetilde{c},
$$
for any $Y \in {\bf C}(n,2r)$ and any chosen $t$.

Hence, $V$ cannot be a tangent vector.

To prove \eqref{assume two}, we also chose $X_{t}$ to be $
X_{t}=\begin{pmatrix}
M+tA & tB\\
-t B^{T} & tD-t^2 B^{T}(M+tA)^{-1}B
\end{pmatrix}$, where ${\rm rank} D>2r-2k$. It can be seen that $X_{t}\notin {\bf C}(n,2r)$ by the Gauss elimination: 
\begin{equation}\label{svd}
X_{t}\sim N_{t}=\begin{pmatrix}
M+tA & 0\\
0 & tD
\end{pmatrix}.
\end{equation}

Then, it suffices to prove 
\begin{equation}\label{assume three}
N_{t} \ is \ not \ at \ distance \ o(t) \ from \ {\bf C}(n,2r),
\end{equation}
since by the norm inequality, when $t$ is small enough, $X_{t}$ is not at distance $o(t)$ from ${\bf C}(n,2r)$.

The proof of \eqref{assume three} depends on the Low Rank Approximation and Eigenvalue Perturbation theories in matrix analysis, readers are referred to \cite{EY36},\cite{HJ13},\cite{SS90} for more informations.

From \eqref{svd}, 
$$
N_{t}N_{t}^{T}=\begin{pmatrix}MM^T+t(AM^{T}+MA^{T})+t^2AA^{T} & 0 \\ 0 & t^2DD^{T} \end{pmatrix}.
$$

Denote the nonascending eigenvalues of the upper left block by: $\mu_{1}\geq \cdots \geq \mu_{2k}$ and the singular values of $D$ by: $\lambda_{1}=\lambda_{1} \geq \cdots \geq \lambda_{s}=\lambda_{s}>0$, where ${\rm rank} \ D:=s>r-k$. Additionally, denote the maximal singular value of $AM^{T}+MA^{T}+tAA^{T}$ by $\sigma$. 

Then from Weyl's theorem (cf. Corollary 4.3.15 in \cite{HJ13}), we can see that for every $1\leq i\leq k$, both
$$
\mu_{2i-1},\mu_{2i}\in [x_{i}^2-t\sigma, x_{i}^2+t\sigma].
$$

It follows that there exists some $t_{0}$, when $t\in (0,t_{0})$, the singular values of $N_{t}$ in \eqref{svd} satisfy:
\begin{equation}\label{assume four}
\sqrt{\mu_{1}}\geq \cdots \geq \sqrt{\mu_{2k}}>t \lambda_{1}=t \lambda_{1} \geq \cdots \geq t \lambda_{s}=t \lambda_{s}>0,
\end{equation}
where $s>r-k$.

Unlike a direct application of Low Rank Approximation for determinantal varieties in \cite{KL99}, we should give a singular value decomposition $(\textit{SVD})$ for skew-symmetric matrices, which can also be achieved from its canonical forms. 

Any nonzero matrix $X\in {\bf C}(n,2r)$ can be written in the canonical form: $X=AM_{0}A^{T}$, $A\in SO(n)$, $M_{0}=\sum_{i=1}^{k}x_{i}X_{2i-1,2i}$ with $x_{1}\geq \cdots \geq x_{k}>0(1\leq k \leq r)$.

It is easy to see every diagonal block $\begin{pmatrix} 0 & x_{i} \\ -x_{i} & 0 \end{pmatrix}$ 
has the following \textit{SVD}(it is not unique for left and right singular vectors!):
$$
\begin{pmatrix} 0 & x_{i} \\ -x_{i} & 0 \end{pmatrix}=\begin{pmatrix}0 & 1 \\ -1 & 0 \end{pmatrix}\begin{pmatrix} x_{i} & 0 \\ 0 & x_{i} \end{pmatrix}\begin{pmatrix} 1 & 0 \\ 0 & 1 \end{pmatrix},
$$
let 
$$
P={\rm diag} \left\lbrace \begin{pmatrix}0 & 1 \\ -1 & 0 \end{pmatrix},\cdots,\begin{pmatrix}0 & 1 \\ -1 & 0 \end{pmatrix},0,\cdots,0 \right\rbrace
$$ where the nonzero diagonal blocks have number $k$.

Then,
$$
X=(AP)\Sigma A^{T}, \Sigma={\rm diag}\{x_{1}\textbf{I}_{2},\ldots,x_{k}\textbf{I}_{2},0,\ldots,0\},
$$
gives the \textit{SVD}.

Now, we can apply{ Low Rank Approximation, combined with \eqref{assume four}, note that $|X|^2=\frac{1}{2} tr(XX^{T})$,we deduce that
$$
{\rm Dist}(N_{t},{\bf C}(n,2r))\geq t \sqrt{(\lambda_{r-k+1}^2+\cdots+\lambda_{s}^2)},
$$
thus proves \eqref{assume three} and finishes all the proofs. $\Box$

\medskip\noindent

Assume $N$ is a nearly regular point of ${\bf C}(2n,2n-2)$, then by Theorem 7.2,
$$
{\rm Tan}({\bf C}(2n,2n-2),N)\cong {\bf C}(4,2) \times \textbf{R}^{(n-2)(2n+3)}.
$$

Recall that ${\bf C}(4,2)$ is non-minimizing. Also, any surface is area-minimizing iff its Cartesian product with any Euclidean space is area-minimizing. It concludes the following

\medskip\noindent

{\bf Theorem 7.4}.  {\it For nearly regular points of Pfaffian hypersurfaces ${\bf C}(2n,2n-2)$, their tangent cones are not area-minimizing. }

\medskip\noindent
\medskip\noindent

{\it Acknowledgments.}\,\,  This project is supported by the NSFC (No.11871445, No.11871450), the project of Stable Support for Youth Team in Basic Research Field, CAS(YSBR-001) and the Fundamental Research Funds for the Central Universities.

\begin{flushleft}
\medskip\noindent
\begin{tabbing}
XXXXXXXXXXXXXXXXXXXXXXXXXX*\=\kill
Hongbin Cui\\
School of Mathematical Sciences, University of Science and Technology of China\\
Wu Wen-Tsun Key Laboratory of Mathematics, USTC, Chinese Academy of Sciences\\
96 Jinzhai Road, Hefei, 230026, Anhui Province, China\\

E-mail: cuihongbin@ustc.edu.cn

\end{tabbing}

\end{flushleft}

\begin{flushleft}
\medskip\noindent
\begin{tabbing}
XXXXXXXXXXXXXXXXXXXXXXXXXX*\=\kill
Xiaoxiang Jiao\\
School of Mathematical Sciences, University of Chinese Academy of Sciences\\
19A Yuquan Road, Beijing, 100049, China\\

E-mail: xxjiao@ucas.ac.cn

\end{tabbing}

\end{flushleft}

\begin{flushleft}
\medskip\noindent
\begin{tabbing}
XXXXXXXXXXXXXXXXXXXXXXXXXX*\=\kill
Xiaowei Xu\\
School of Mathematical Sciences, University of Science and Technology of China\\
Wu Wen-Tsun Key Laboratory of Mathematics, USTC, Chinese Academy of Sciences \\
96 Jinzhai Road, Hefei, 230026, Anhui Province, China\\
E-mail: xwxu09@ustc.edu.cn

\end{tabbing}

\end{flushleft}

\end{document}